\documentclass[a4,11pt]{amsart}
\usepackage[arrow,matrix]{xy}
\usepackage{amssymb,yfonts}
\usepackage{amsbsy,amsmath}
\newtheorem{thm}{Theorem}[section]
\newtheorem{lemma}[thm]{Lemma}
\newtheorem{def+lem}[thm]{Definition+Lemma}
\newtheorem{cor}[thm]{Corollary}
\newtheorem{con}[thm]{Construction}
\newtheorem{definition}[thm]{Definition}
\newtheorem{remark}[thm]{Remark}
\newcommand{\re}{\mbox{\rm Re}}
\newcommand{\im}{\mbox{\rm Im}}

\newcommand{\dd}{\partial\bar\partial}
\newcommand{\seq}{\longrightarrow}

\newcommand{\imp}{\Rightarrow}

\newcommand{\rk}{\mbox{rk}}
\newcommand{\bb}{\mathbb}

\newcommand{\codim}{\mbox{\rm codim }}
\newcommand{\Aut}{\mbox{\rm Aut}}
\newcommand{\Alb}{\mbox{\rm Alb}}
\newtheorem{example}[thm]{Example}
\newenvironment{impl}{\vspace{3pt}\indent
                       \textsc{\it Implementation.}\quad }
                       {\hfill$\square$\vspace{3pt}}
\newenvironment{rem}{\begin{remark}\rm}{\end{remark}}
\newenvironment{defn}{\begin{definition}}{\end{definition}}

\title[Invariance and Ricci-flatness of metrics on open manifolds]{On invariance and Ricci-flatness of Hermitian metrics on open manifolds}

\author{Bert K\"ohler and Marco K\"uhnel}
\address{Bert K\"ohler\\Philipps-Universit\"at Marburg\\Department of Mathematics\\35032 Marburg\\Germany}
\email{koehlerb@mathematik.uni-marburg.de}
\address{Marco K\"uhnel\\Otto-von-Guericke-Universit\"at Magdeburg\\FMA / IAN\\Postfach 4120\\
39016 Magdeburg\\Germany}
\address{Current address of M. K.: Mathematics Institute\\University of Warwick\\Coventry CV4 7AL\\United Kingdom}
\email{mkuehnel@maths.warwick.ac.uk}
\date{\today}
\thanks{The authors acknowledge gratefully support by the 
DFG priority program 'Global Methods in Complex Geometry'.}

\begin{document}

\begin{abstract}We discuss a technique to construct Ricci-flat hermitian metrics on complements of
(some) anticanonical divisors of almost homogeneous manifolds and inquire into when this metric is complete
and K\"ahler.
This construction has a strong interplay with invariance groups of 
the same dimension as the manifold acting with an open orbit. Lie groups of this type we call divisorial.
As an application we describe compact manifolds admitting a 
divisorially invariant
K\"ahler metric on an open subset. Finally, we see a connection between the reducibility
of the anticanonical divisor and the non-triviality of the K\"ahler cone on the complement.
\end{abstract}

\maketitle

\setcounter{section}{-1}
\section{Introduction}

The first problem addressed in this paper is the construction of Ricci-flat metrics on open manifolds.
As a model for this situation serves the complement of a divisor on a compact manifold $X$. If
$D\in|-K_X|$, then the section of $-K_X$ vanishing exactly on $D$ yields an isomorphism 
$\Omega^n_{X\setminus D}\cong{\mathcal O}_{X\setminus D}$. In analogy to the Calabi conjecture on compact 
manifolds this raises
the expectation that there exists a complete Ricci-flat K\"ahler metric on $X\setminus D$. Moreover, methods
to find such a metric may also work, if $D$ is not reduced, leading to the speculative existence of
Ricci-flat K\"ahler metrics on $X\setminus D$, whenever $-(D+K_X)$ is effective. 
Part of this program has already been established. Tian and Yau have proved in \cite{TY} and \cite{TY2}
the existence of a complete Ricci-flat K\"ahler metric in case $D\in|-K_X|$ is neat, almost ample and smooth.
Bando and Kobayashi have shown the claim, if $rD\in|-K_X|$ for $r>1$ and $D$ is ample, smooth and admits
on itself a K\"ahler-Einstein metric. The metrics involved in the construction contain logarithmic terms.  
Of course, the techniques introduced in \cite{TY}, \cite{TY2} and \cite{BK} cannot be easily 
generalized to the reducible or non-reduced case. The method described in Section 1 does not care
about this and is compatible with an algebraic structure. In particular, the metrics involved can be described
in terms of polynomials, if $X$ is algebraic. Of course, this nice structure has a high prize: either
$X\setminus D$ has trivial geometry and $D$ is expected to be 'very' reducible, 
or we have to drop the K\"ahler condition. Nevertheless, this approach
shows some fundamental differences between the reducible and the smooth case of the divisor $D$.

The other problem is determining highly symmetric metrics on $X\setminus D$. Both problems get related by
the idea that Ricci-flatness should be forced by a high order of symmetry. For example, a naive calculation
shows that all metrics on ${\bb P}^2\setminus\{\mbox{3 general lines}\}$ invariant under the
connected invariance group of the three lines are complete, Ricci-flat and K\"ahler. 
In Chapter 2 we
explain the connection between the construction in Chapter 1 and the symmetries of $D$ resp. symmetries of
the metric. Here we discuss continuous symmetries. A striking point is that the
metric is K\"ahler if and only if it is symmetric. Moreover, in this case the symmetry group is abelian.
This allows a description of the manifolds admitting a divisorially invariant K\"ahler metric
on an open subset; in particular, we recognize $D$ to be reducible or non-reduced, if $X$ is homogeneous and projective. 
By divisorial invariance we mean that the action of $G$ has an open orbit and
$\dim G=\dim X$, if $G$ denotes the symmetry group of 
$D$. Three general lines in ${\bb P}^2$ satisfy this condition. Parts of the description are well known.

Winkelmann treated in \cite{wi} the problem when $T_X(-\log D)$ is trivial. Of course, the condition that there is a Ricci-flat metric on $X\setminus D$ is much weaker
than the triviality of $T_X(-\log D)$. However, in the respective K\"ahler cases there are great similarities. We will note this at the appropriate place.
 
As a last topic we inquire into a sort of K\"ahler classes of the constructed metrics. 
Two metrics shall be regarded as equivalent if they differ only by a K\"ahler potential. 
We call the cone of $G$-invariant metrics generated 
by this equivalence $K_G(X,D)$. If we denote $n:=\dim X$, then we will prove that $\dim K_G(X,D)\ge\frac 12n(n-1)$. So even if $\dim \Alb(X)=0$ the cone $K_G(X,D)$
is highly non-trivial. This effect for $K(X,D)$ 
is in close relation with the reducibility of $D$. In the appendix we will show that $K(X,D)=0$, if
$D$ is smooth and ample and $X$ has simple enough topology, e.g. $X={\bb P}^3$.

\section{Vocabulary}

We consider compact complex manifolds $X$. Of great importance will be the automorphism group $\Aut(X)$
and its action on $X$. If $G\subset \Aut(X)$ is a Lie group, we write $G^0$ for the connected component of 
$G$ containing
the identity. By $\textswab g$ we denote the Lie algebra of $G$.
If $D\subset X$ and $g$ is a metric on $X\setminus D$, then we define
$$\Aut(X,D):=\{\phi\in \Aut(X)|\,\,\phi|_D\in \Aut(D)\}$$
and
$$\Aut(X,g):=\{\phi\in \Aut(X,D)|\phi^*g=g\}.$$
If $Y$ is some complex manifold and $g$ a metric on $Y$, we also denote
$$\Aut(Y,g):=\{\phi\in \Aut(Y)|\phi^*g=g\}.$$

In most cases we will further assume that $X$ is
{\em almost homogeneous}. We will often make use of the following equivalences.

\begin{def+lem}A compact complex manifold $X$ is called {\em almost homogeneous}, if there is a
Lie group $G\subset \Aut(X)$ such that one (and then all)
of the following properties are satisfied:
\begin{enumerate}
\item The action of $G$ has an open orbit,
\item the action of $G^0$ has an open orbit,
\item ${\textswab g}:=T_1G$ generates $T_X$ at the general point,
\item there is a vector space $V\subset\textswab g$ with $\dim V=\dim X$, which generates $T_X$ at
the general point.
\end{enumerate}
\end{def+lem}

If $G\subset \Aut^0(X)$ is a Lie group which has an open orbit, then we say $G$ {\em acts almost
transitively} on $X$.

In the K\"ahler case we will encounter a special form of abelian Lie groups, so called {\em semi-tori}. 

\begin{defn}A complex Lie group $G$ is a {\em semi-torus}, if there is a number $n$ and a discrete subgroup $\Lambda\subset{\bb C}^n$ such that
$\Lambda\otimes{\bb C}={\bb C}^n$ and $G={\bb C}^n/\Lambda$.\end{defn}

The other aspect of the paper is K\"ahler metrics on complements of divisors. If $Y$ is any
complex manifold and $G\subset \Aut(Y)$
is a Lie group, we denote by $M_G(Y)$ the set of all K\"ahler metrics on $Y$ invariant
under $G$. If $X$ is a compact complex manifold, $D\subset X$ a divisor, $Y=X\setminus D$ and
$G\subset \Aut(X,D)$, then we also write $M_G(X,D):=M_G(Y)$ to emphasize that we use only automorphisms
of $X$ for the construction.

For analytical application a K\"ahler potential is very useful rather than a 
$\overline\partial$-exact $(1,1)$-form. In the compact case there is no difference, but in the open case
we have to distinguish. This is the reason why we do not consider the K\"ahler cone as a subcone
of $H^{1,1}(X)$, but use a finer equivalence relation. 

\begin{defn}For K\"ahler metrics $g,g'$ on a complex manifold $Y$
$$g\sim g':\iff \,\exists\phi\in C^\infty(Y,{\bb R}): \omega_g-\omega_{g'}=i\dd\phi$$
is an equivalence relation. If $G\subset \Aut(Y)$, we define
$$K_G(Y):=M_G(Y)/\sim$$
and call it the {\em $G$-K\"ahler cone} of $X\setminus D$. We also abbreviate
$K(Y):=K_1(Y)$ and call this cone the {\em K\"ahler cone} of $Y$. Note that this definition 
of the K\"ahler cone coincides with the usual, if $Y$ is compact.
Like above, if $X$ is a compact
complex manifold, $D\subset X$ a divisor, $Y=X\setminus D$ and $G\subset \Aut(X,D)$, we also write
$K_G(X,D):=K_G(Y)$. 
\end{defn}

Since parts of the paper are written in the language of differential geometry, we use
co- and contravariant indexing conventions as well as Einstein's sum convention. 
In order to be able to do this, we distinguish
indices arising from non-differential context by setting them in brackets, if appropriate. 
The decision, whether such an index is sub- or superscript, is made by considering the beauty of 
involved formulas.
For example, a set of vector fields is denoted by $s^{(i)}$. The same vector fields in local 
coordinates will be written $s^{ik}\frac{\partial}{\partial z^k}$. Here we omit the brackets, because
we want to put the components into a matrix. The only unlucky point of this convention is where powers
of coordinates appear. But we believe that also in these cases the meaning will become clear.

\section{Metrics generated by vector fields}

\subsection{Construction of the metrics}

In this section we discuss a method to construct complete, Ricci-flat hermitian metrics. 
These are neither necessarily invariant nor K\"ahler. We will determine the cases when the metric is
K\"ahler.

The construction of the divisor is widely used in many works about almost homogeneous manifolds. However, in this place we concentrate on the metric which comes
with the construction. For this reason we give here a detailed description. 

We use the notion of an abelian subspace of $H^0(T_X)$. We call $V\subset H^0(T_X)$ abelian, if for all $\zeta,\xi\in V$ holds $[\zeta,\xi]=0$. In general, however, we do not require $V$ to be an algebra.

\begin{con}Let $X$ be an almost homogeneous compact complex manifold. Let
${\mathcal B}=\{s^{(1)},...,s^{(n)}\}\subset H^0(T_X)$ generate $T_X$ in
the general point and denote $V=<{\mathcal B}>$ the vector space spanned by ${\mathcal B}$. This yields a divisor 
$D_V\in|-K_X|$ and a Ricci-flat hermitian metric $g_{\mathcal B}$ on $X\setminus D_V$. The metric $g_{\mathcal B}$
is K\"ahler if and only if $V$ is abelian.
\end{con} 

\begin{impl}Since $s^{(1)},...,s^{(n)}\in H^0(T_X)$ generate $T_X$ in a general point,
$\bigwedge_{i=1}^ns^{(i)}$ vanishes exactly on a divisor $D_V$. Obviously $D_V\in|-K_X|$. Since on $X\setminus D_V$
the $s^{(i)}(x)$ form a basis of $T_{X,x}$, we may construct $s_{(i)}\in\overline{T}_{X,x}^*$ by
prescribing $s_{(i)}(\overline{s^{(j)}})=\delta_{ij}$ on $X\setminus D_V$. Further we can extend this correspondence to a linear map
$${}^\dagger:T_X\seq\overline{T}_X^*$$  and
define 
$$g_{\mathcal B}(s\otimes t):=s^\dagger(t),$$
if $s\in T_{X,x}, t\in \overline{T}_{X,x}$. In a local chart we denote
$s^{(i)}=s^{ik}\frac{\partial}{\partial z^k}.$ We denote by $(s_{ij})=\sigma$ the inverse matrix of $(s^{ij})$.
Then
$$g_{{\mathcal B},ij}=g_{\mathcal B}\left(\frac{\partial}{\partial z^i}\otimes\frac{\partial}{\partial \overline{z^j}}\right)=s_{ik}\overline{s_{jk}},$$
what yields
$$Ric(g_{\mathcal B})=\frac{i}{2\pi}\partial\overline\partial\log\det g_{\mathcal B}=\frac{i}{2\pi}\partial\overline\partial\log\det\sigma+
\frac{i}{2\pi}\partial\overline\partial\log\overline{\det\sigma}=0,$$
since $\det\sigma$ is holomorphic.

A short calculation shows that $g_{\mathcal B}$ is K\"ahler if and only if
$$s_{ij,l}=s_{lj,i}$$
for all $i,j,l$. Converting this condition to the vector field components yields
$$s^{ij}s^{kl}_{,j}=s^{kj}s^{il}_{,j}$$
for all $i,k,l$. Of course, this is the condition
$$[s^{(i)},s^{(k)}]=0$$
for all $i,k$. But this means that $V$ is abelian. 
\end{impl}

Note that by construction $T_{X\setminus D_V}$ is trivial. In \cite{wi} the problem is addressed when $T_X(-\log D)$ is trivial and answered in terms of
the existence and action of a semi-torus. In the next section we will be able to describe this property in terms of $g_{\mathcal B}$ and $D_V$. 

\subsection{Completeness of the metrics}

Since we are dealing with open manifolds we should address the problem of completeness of the constructed
metric. For this purpose we introduce some new notation. First, we define $S:=(s^{ik})=\sigma^{-1}$. Then
we interpret $S^t:{\mathcal O}_X^{\oplus n}\seq T_X$ as a sheaf homomorphism
and define 
$${\mathcal L}:=\ker S^t.$$
Of course, ${\mathcal L}$ is supported on $D_V$ and we prove

\begin{lemma} If $D_V$ is smooth, then ${\mathcal L}$ is a line bundle on $D_V$.\end{lemma}

\begin{proof}${\mathcal L}$ is line bundle if and only if $\rk S|D_V=n-1$ everywhere.  So assume, that
in $x\in D_V$ we have $\rk S(x)<n-1$. Then all $(n-1)\times(n-1)$-minors of $S$ vanish in $x$,
in particular $d\det S(x)=0$, hence $x\in Sing(D_V)$.
\end{proof}

Now assume that $D_V$ is smooth and 
consider $\omega_{(i)}:=\bigwedge_{j\not= i}s^{(j)}|_{D_V}\in H^0(\bigwedge^{n-1}T_X|D_V)=H^0(\Omega^1_X|D_V
\otimes N_{D_V|X})$. These are related via ${\mathcal L}$ by the equation
$$(-1)^j\lambda_{(j)}\omega_{(i)}+(-1)^i\lambda_{(i)}\omega_{(j)}=0$$
for all $i,j$ and $x\in D_V,\lambda\in{\mathcal L}_x$. Again, smoothness of $D_V$ implies that in
a point $x\in D_V$ not all $\omega_{(i)}$ can vanish. 
Hence the vector spaces ${\mathcal F}_x:=<(\omega_{(i)})_i>$ 
are one-dimensional and form a line bundle ${\mathcal F}\subset\Omega^1_X|D_V\otimes N_{D_V|X}$. By
looking at the natural local trivializations of ${\mathcal L}$ and ${\mathcal F}$ it is easy to see that
$${\mathcal F}\cong{\mathcal L}^\vee.$$

Note that the inclusion $i:D_V\seq X$ yields via Poincar\'e Duality a homomorphism $i_*:H^*(D_V,{\bb R})\seq
H^*(X,{\bb R})$ of degree $2$. With this notation in mind, the very definition of ${\mathcal F}$ implies that
$$i_*c_1({\mathcal F})=c_2(X).$$

The tensored dual tangent sequence
$$0\seq{\mathcal O}_{D_V}\seq\Omega^1_X|D_V\otimes N_{D_V|X}\stackrel{\pi}{\seq}\Omega^1_{D_V}\otimes N_{D_V|X}\seq 0,$$ 
allows us to formulate the property ${\mathcal F}={\mathcal O}_{D_V}=\ker\pi$.

\begin{lemma}\label{ker}If $D_V$ is smooth, then $g_{\mathcal B}$ is complete if and only if ${\mathcal F}=\ker\pi$.\end{lemma}

\begin{proof}Let us choose $x\in D_V$ and local coordinates in a small open subset $U\subset X$ 
such that $D_V=\{z^1=0\}$. Furthermore, denote $U':=U\cap D_V$ and $pr:U\seq U'$ the projection
induced by the local coordinates.
choose $0\not=\lambda\in{\mathcal L}(U')$ and an order of ${\mathcal B}$ such that $\lambda_{(1)}\equiv 1$.
If we now define $B$ by 
$$B_{ij}(z):=\left\{\begin{array}{rl}1&\mbox{ if }i=j\\
-\lambda_{(j)}(pr(z))&\mbox{ if } i=1,j\not=1\\
0&\mbox{ else } \end{array}\right.,$$
then $\tilde S:=BS$ is just $S$ replaced by a first row vanishing on $U'$.
Now we have a look at $\tilde\sigma:=(\tilde S)^{-1}$. Since the first row is
identically $0$ on $U'$, we conclude that $\tilde s_{1i}\in{\mathcal O}(U)$ for $i>1$. Since
$$g_{11}=\sum_{i>1}|\tilde s_{1i}|^2-2\re(\overline{\tilde s_{11}}\sum_{i>1}\lambda_{(i)}
\tilde s_{1i})+
|\tilde s_{11}|^2\sum|\lambda_{(j)}|^2,$$
we see now, that $g$ is complete if and only if $\tilde s_{11}\sim\frac{1}{z^1}$ for all such choices of
coordinates. Indeed, $\tilde s_{11}=s_{11}$. Hence,
if we denote $A_{ij}$ the $(i,j)$-entry of the cofactor matrix of $S$, then the condition is
equivalent to $A_{11}(x)\not=0$. If we now choose other coordinates $z^{\prime 1},...,z^{\prime n}$
such that $D_V=\{z^{\prime 1}=0\}$ and denote $J:=(\frac{\partial z^{\prime j}}{\partial z^i})_{ij}$, 
$h:=\frac{z^{\prime 1}}{z^1}\not= 0$, then
$$A'_{11}=\det J(h^{-1}A_{11}+\sum_{k>1}J^{-1}_{1k}A_{1k})\not= 0$$
Since the coordinate transform was arbitrary, we conclude that
$$A_{1k}(x)=0\mbox{ for all }k>1.$$
This is equivalent to
$0\not=\omega_{(1)}\in H^0(\ker\pi).$
This again means ${\mathcal F}=\ker\pi$.
\end{proof} 

Note that completeness as well as the K\"ahler property of $g_{\mathcal B}$ depend only on $V$.

If $D_V$ is not smooth, then we consider $D_V^0$, the regular part of $D_V$ and the corresponding
objects ${\mathcal F}^0, {\mathcal L}^0, \pi^0$, which are obtained by restriction to $D_V^0$. We now show
the preceding Lemma for the singular case.

\begin{lemma}$g_{\mathcal B}$ is complete if and only if ${\mathcal F}^0=\ker\pi^0$.\end{lemma}

\begin{proof}'$\Rightarrow$': If $g_{\mathcal B}$ is complete, the same arguments as in Lemma \ref{ker}
imply that ${\mathcal F}^0=\ker\pi^0$.

'$\Leftarrow$': Let locally $D_V=\{f=0\}$ in a small open neighborhood $U\subset X$.
Then we can choose functions $z^2,...z^n$ which give local coordinates together with $f$ on
the set $\tilde U:=U\setminus\{df\wedge dz^2\wedge...\wedge dz^n=0\}$. Like in the proof above
we argue that $g_{\mathcal B}$ is complete if $A_{11}(x)\not=0$ for choices like above and $x\in D_V$.
Since by assumption this is true for $x\in D_V\setminus Sing(D_V)$, extension of the
holomorphic function $A_{11}$ to $U'$ yields a non-zero function $A_{11}\in{\mathcal O}^*(U')$, if
$D_V$ was normal. If $D_V$ is not normal, we choose an embedded normalization
$$\xymatrix{\tilde D_V\ar[r]^{\nu}\ar[d]_j&D_V\ar[d]^i\\
\tilde X\ar[r]_\mu&X},$$
where $i,j$ denote inclusions and $\nu$ the normalization of $D_V$. Now we apply the same arguments to
the pseudometric $\mu^*g$ and obtain by looking at paths $\gamma$ such that $\mu|_\gamma$ is a
diffeomorphism that $g$ is complete.  
\end{proof}

Now we can see a connection to the invariance group. 

\begin{lemma}\label{c1}$g_{\mathcal B}$ is complete if and only if 
$\exp(V)\subset \Aut^0(X,D_V)$.
\end{lemma}

\begin{proof}
If we denote by ${\mathcal V}$ the sheaf on $D$ generated by ${\mathcal V}_x:=<\{s^{(i)}(x)\}>$ for $x\in D$,
then ${\mathcal F}({\mathcal V})\equiv 0$, if we regard ${\mathcal F}\subset Hom(T_X|_{D_V},N_{{D_V}|X})$. 
Since ${\mathcal F}^0=\ker\pi^0$, figuring out the dualized maps
$$0\seq T_{D_V^0}\seq T_X|D_V^0\stackrel{p^0}{\seq}{\mathcal F}^{0\vee}\otimes N_{{D_V^0}|X}\seq 0$$
yields ${\mathcal V}^0=\ker p^0=T_{D_V^0}.$ This is equivalent to $V|_{D_V^0}\subset H^0(T_{D_V^0})$. This 
again means that every
$\phi\in\exp(V)$ holds $D_V$ invariant, so $\exp(V)\subset \Aut^0(X,D_V)$ is an equivalent condition.
\end{proof}

Taking into account the results of Section \ref{sym}, which are obtained independently of the considerations
about completeness, we even find

\begin{thm}\label{complete}$g_{\mathcal B}$ is complete if and only if $V$ is a Lie subalgebra.
\end{thm}

\begin{proof} 
'$\Rightarrow$':  Since $g_{\mathcal B}$ is complete, by Lemma \ref{c1}  we obtain $\exp(V)\subset \Aut^0(X,D_V)$.
Hence $\Aut^0(X,D_V)$
acts almost transitively. So Lemma \ref{max} yields $\dim\exp(V)=\dim X=\dim \Aut^0(X,D_V)$. This means
$$T_1 \Aut^0(X,D_V)=T_1\exp(V)=V,$$ hence $V$ is a Lie subalgebra. 

'$\Leftarrow$': If $V$ is a Lie subalgebra, Lemma \ref{prebinding} implies the desired property  $\exp(V)\subset \Aut^0(X,D_V)$.
\end{proof}

Now it is clear from the previous arguments that $T_X(-\log D_V)$ is trivial, if $g_{\mathcal B}$ is complete and $D_V$ is a simple normal crossings divisor. In Chapter
\ref{kaehler} we will show furthermore that $G:=\exp(V)$ is a semi-torus, if $g_{\mathcal B}$ is complete and $D_V$ is reduced.
 
The existence of a complete $g_{\mathcal B}$ for a smooth $D_V$ restricts the geometry of $X$
significantly:

\begin{cor}If $D_V$ is smooth and $g_{\mathcal B}$ is complete, then $c_2(X)=0$.\end{cor}

\begin{proof}
${\mathcal F}=\ker\pi\cong{\mathcal O}_{D_V}$, hence $c_2(X)=i_*c_1({\mathcal O}_{D_V})=0.$
\end{proof} 

Now we have seen 
that $c_2(X)\not=0$ implies that the divisor $D_V$ is singular, if $V$ is a Lie subalgebra. We will
see later that the K\"ahler condition allows an explicit description of the singularities, at least on projective homogeneous manifolds.

\section{Symmetries of the divisor and the metric}\label{sym}

In this section we want to relate the construction above to the appearance of symmetries on $D$ and the
metric.
As it may be not hard to guess, this connection is made by Lie theory. Throughout 
this section we mean always a complex Lie group when we speak of Lie groups. 

If $G$ is a Lie group we identify $\textswab g=T_1G$, and if $G\subset \Aut^0(X)$, then we furthermore identify
$\textswab g$ with the subvector space of $H^0(T_X)$ given by the vector fields 
$s(x):=\frac{\partial}{\partial t}g(t)x|_{t=0}$, where $g(t)$ denotes a (holomorphic) 
path in $G$ with $g(0)=1$ and
$\frac{\partial}{\partial t}g(t)|_{t=0}=\xi\in T_1G$. Furthermore we have an action of $G$ on $T_1G$
by $h\xi:=\frac{\partial}{\partial t}hg(t)h^{-1}|_{t=0}$, if $h\in G$.

\begin{lemma}\label{prebinding}
Let $X$ be a compact complex manifold of dimension $n$, $G\subset \Aut^0(X)$ a Lie group 
acting almost transitively on $X$ and $\textswab g$ be the corresponding Lie algebra. Then

\begin{enumerate}
\item $D\in|-K_X|$ and $G\subset \Aut^0(X,D)\imp\,\, D=D_V$ for all $V\subset{\textswab g}$ with $\dim V=n$ and
generating $T_X$ in the general point,
\item if $\dim G=n$, then $G\subset \Aut^0(X,D_{\textswab g})$,
\item if $\dim G=n$, ${\mathcal B}\subset\textswab g$ is a basis, then:
$G\subset \Aut^0(X,g_{\mathcal B})\iff G\mbox{ is abelian.}$
\end{enumerate}
\end{lemma}

\begin{proof}'(i)':Let $V\subset\textswab g$ be an $n$-dimensional vector space generating $T_X$ in the
general point and ${\mathcal B}\subset V$ a basis. Since $D$ is $G$-invariant, for 
any $s\in\textswab g\subset H^0(T_X)$ the restriction $s|_D$ gives an element of $H^0(T_D)$. 
Since $\dim D=n-1$, this implies $\bigwedge_{s\in{\mathcal B}}s|_D=0$, hence $D\subset D_V$. But $D$ and $D'$
are both elements of $|-K_X|$, hence $D=D_V$.

'(ii)': If $s\in\textswab g\subset H^0(T_X)$ is given by $\xi\in T_1G$, then for $h\in G$ the pullback
$h^*s$ is given by $h^{-1}\xi\in T_1G$, hence $h^*s\in\textswab g$. Furthermore $h^*$ maps a basis of
$\textswab g$ to a basis again, because $h$ is an automorphism.
This proves that
$\bigwedge s^{(i)}=0\iff \bigwedge h^*s^{(i)}=0,$
if $s^{(1)},...,s^{(n)}$ is a basis of $\textswab g$. 
Hence $h(D)=D$. This proves that $D$ is $G$-invariant.

'(iii)': Let $\omega$ be the fundamental $(1,1)$-form of $g_{\mathcal B}$.
Of course, $g_{\mathcal B}$ is $G$-invariant if and only if ${\mathcal L}_s\omega=0$, if $s\in\textswab g$ and 
${\mathcal L}$ denotes the Lie derivative. If $C_s$ denotes the contraction by $s$, then 
${\mathcal L}_s=dC_s+C_sd$ (see e.g. \cite[V,5]{lang}). Let us choose local coordinates like in the construction.
Then
$$0={\mathcal L}_{s^{(l)}}\omega=\sum_{i,j,k,m}s^{lm}(s_{ik,m}-s_{mk,i})\overline{s}_{jk}dz^i\wedge 
d\overline{z}^j.$$
Since $S$ is invertible on $X\setminus D$,
we obtain
$$s_{ik,m}-s_{mk,i}=0$$
for all $i,k,m$.
This we identified at an earlier place as the condition $$[s^{(i)},s^{(j)}]=0$$ for all $i,j$.
Hence $\textswab g$ is abelian. It is well known that 
this is equivalent to $G$ to be abelian.
\end{proof}

As a first application, we obtain kind of uniqueness of the vector space $V$ in the construction.

\begin{lemma}\label{max}Let $X$ be a compact complex manifold and $D\in|-K_X|$. 
If $\Aut(X,D)$ acts almost transitively, then $\dim \Aut(X,D)=\dim X$. 
\end{lemma}

\begin{proof}We abbreviate $G:=\Aut(X,D)$ and $n:=\dim X$.
Of course, $\dim G\ge n$, since otherwise $T_1G$ could not generate $T_X$ in any point. 
Now choose
${\mathcal B}:=\{s^{(1)},...,s^{(n+1)}\}\subset{\textswab g}$ such that ${\mathcal B}\setminus\{s^{(n+1)}\}$ generates
$T_X$ in the general point. Further denote $\eta_{(i)}:=\bigwedge_{j\not=i}s^{(j)}$.
and
by $V_{(i)}$ the vector space generated by ${\mathcal B}\setminus\{s^{(i)}\}$. 

Since $\dim X=n$, we can find meromorphic functions $f_{(i)}\in{\mathcal M}_X(X)$ such that
$s^{(n+1)}=\sum_{i=1}^nf_{(i)}s^{(i)}$. By assumption we know $\eta_{(n+1)}\not\equiv 0$. 

If $\eta_{(i)}\equiv 0$,
then we see by $\eta_{(i)}=f_{(i)}\eta_{(n+1)}$, that $f_{(i)}\equiv 0$. 

If $\eta_{(i)}\not\equiv 0$, then $V_{(i)}$ generates $T_X$ in the general point 
and we may use Lemma \ref{prebinding} to obtain
$$\eta_{(n+1)}=0\iff\eta_{(i)}=0\iff f_{(i)}\eta_{(n+1)}=0.$$
Hence $f_{(i)}$ has no zeroes. By exchanging $s^{(i)}$ and $s^{(n+1)}$ we also see that
$f_{(i)}$ has no poles. Hence $f_{(i)}$ is constant.

Now we proved that every $f_i$ is constant, hence
$s^{(1)},...,s^{(n+1)}$ are linearly dependent and $\dim G=n$.
\end{proof}

Note that the connection between the invariance group and the anticanonical system is essential.
For example, the invariance group of a point in ${\bf P}^1$ is two-dimensional and acts almost
transitively.

Since every $\dim X$-dimensional Lie group which acts almost transitively yields an invariant
$D\in|-K_X|$, Lemma \ref{max} suggests the following definition.

\begin{defn}If $X$ is compact complex manifold and $G\subset \Aut(X)$ a Lie group with
$\dim G=\dim X$ acting almost transitively on $X$, we say $G$ is a {\em divisorial}
group. If on the other hand $D\in|-K_X|$, we say $D$ has a {\em divisorial invariance
group}, if $\Aut(X,D)$ acts almost transitively. 
Any object invariant under a divisorial group we call {\em divisorially invariant}.
\end{defn} 

Lemma \ref{max} now allows a stronger and more compact formulation of Lemma \ref{prebinding}.

\begin{thm}\label{binding}
Let $X$ be a compact complex manifold, $G$ a divisorial group and $\textswab g$ the corresponding 
Lie algebra. Then for a
divisor $D\in |-K_X|$ holds
$$G=\Aut^0(X,D)\iff\,\, D=D_{\textswab g}.$$
If ${\mathcal B}\subset{\textswab g}$ is a basis, then
$$G=\Aut^0(X,g_{\mathcal B})\iff G \mbox{ is abelian.}$$ 
\end{thm}

\begin{proof}If $D=D_{\textswab g}$, by Lemma \ref{prebinding} $G\subset \Aut^0(X,D)$. Hence $\Aut^0(X,D)$
acts almost transitively and by Lemma
\ref{max} we obtain $\dim G=\dim X=\dim \Aut^0(X,D)$, hence $G=\Aut^0(X,D)$.
\end{proof}

In this context it is appropriate to introduce the notion of a {\em homogeneous pair}.

\begin{defn}A {\em homogeneous pair} $(X,D)$ consists of a compact complex manifold $X$ and
an effective divisor $D$ such that $\Aut^0(X,D)$ acts transitively on $X\setminus D$. 
We call a homogeneous pair $(X,D)$ anticanonical, if $D\in|-K_X|$.
\end{defn}

With this definition we can formulate an equivalence of small categories:

\begin{cor}Given a compact complex manifold $X$, 
there is a one-to-one correspondence between anticanonical homogeneous
pairs of $X$ and divisorial automorphism groups of $X$, given by
$$(X,D)\mapsto \Aut^0(X,D)$$ resp.
$$G\mapsto (X,D_{\textswab g}).$$
\end{cor}

\begin{rem}
\begin{itemize}
\item Note that the proof of Lemma \ref{prebinding} 
also shows that every analytical $G$-invariant set $S$ is contained in
the $G$-invariant $D\in|-K_X|$. 
\item Note that the vector field method is much more general than the invariance approach: There is no
need for the vector space $V\subset H^0(T_X)$ to be an algebra, whereas invariant divisors correspond to
Lie subalgebras of $H^0(T_X)$. 
\item However, if $g$ arises by the general vector field method and is K\"ahler, 
we have $V$ proved to be abelian, in particular $V$ is a Lie subalgebra.
Of course, $\exp:H^0(T_X)\seq \Aut^0(X)$ restricted to $V$ maps to an $n$-dimensional Lie subgroup $G$ leaving
$D$ and $g$ invariant. 
\end{itemize}
\end{rem}

Now we also see that divisorial invariance is exactly the property we had in mind when we expected
that Ricci-flatness should be implied by a high order of symmetry.

\begin{cor}Let $X$ be an $n$-dimensional compact complex manifold, $G\subset \Aut^0(X)$ a divisorial (abelian) 
Lie group. Then there is a
complete Ricci-flat (K\"ahlerian $G$-invariant) metric on $X\setminus D_{\textswab g}$.
\end{cor}

In \cite{wi} Winkelmann proved that $T_X(-\log D)$ is even holomorphically trivial, if $G$ is a complex semi-torus and acts with
only semi-tori as isotropy groups. We will see in the next paragraph that $G$ being a semi-torus is implied by $D_{\textswab g}$
being reduced.

\section{The structure of the K\"ahler case}\label{kaehler}

\subsection{Description of the manifolds}

We saw that the metric $g$ on $X\setminus D$ constructed by an automorphism group $G$ is complete and it is
K\"ahler if 
and only if $G$ is abelian. In this case $g$ is also $G$-invariant. Now we want so see that this construction
yields all $G$-invariant K\"ahler metrics. The first step is to show that invariance groups of K\"ahler 
metrics are abelian. Similar connections between the K\"ahler property and abelian groups are well-known. However, most results in this direction
use compactness of the K\"ahler manifold by employing that all holomorphic one-forms are closed. This formulation makes only use of the K\"ahler form and hence is also
valid in the non-compact case.

\begin{lemma}\label{abel}Let $Y$ be a complex manifold and $g$ a K\"ahler
metric on $Y$. Then $\Aut^0(Y,g)$ is abelian. 
\end{lemma}

\begin{proof}We abbreviate $G:=\Aut^0(Y,g)$. 
Let $\omega$ denote the K\"ahler form of $g$ and $\textswab g\subset H^0(T_Y)$ the Lie algebra 
of $G$. Since $g$ is $G$-invariant, for all $s\in\textswab g$ we obtain $${\mathcal L}_s\omega=0,$$
where ${\mathcal L}$ denotes the Lie derivative. If furthermore $C$ denotes the contraction by
the subscript vector field, ${\mathcal L}_s=dC_s+C_sd$ and hence we conclude $dC_s\omega=0$ for all 
$s\in\textswab g$. Again using an elementary formula (see e.g. \cite[V,5]{lang}) and $d\omega=0$
we obtain for $s,t\in\textswab g$
$$C_{[s,t]}\omega=({\mathcal L}_sC_t-C_t{\mathcal L}_s)\omega={\mathcal L}_sC_t\omega=dC_sC_t\omega+C_sdC_t\omega.$$
Since $\omega$ is a $(1,1)$-form and $s,t$ are holomorphic, $C_sC_t\omega=0$. We already saw that
$dC_t\omega=0$, hence both summands of the right hand side vanish, yielding $C_{[s,t]}\omega=0$.
In local coordinates this means
$$g_{\alpha\bar\beta}[s,t]^\alpha=0.$$
Since the matrix $g_{\alpha\bar\beta}$ is invertible this implies $[s,t]=0$. Hence $\textswab g$
is abelian and therefore also $G$ is abelian.
\end{proof}

It is known (cf, \cite[p. 12]{on1}) that an orbit of $G$ is locally (in $G$) a submanifold. 
A closer look at the argument in the second part of the proof of Theorem \ref{binding} reveals that
any orbit, whose local dimension is smaller than $n$ must be contained in $D_{\textswab g}$. Hence $X\setminus D_{\textswab g}$
is the unique open orbit of $G$. Let $x_0\in X\setminus D_{\textswab g}$ and $\alpha:G\seq X\setminus D_{\textswab g}$ be the action map
$g\mapsto gx_0$. Again it is known (cf. \cite{on1}) that $\alpha$ has constant rank. Since $\alpha$ is
surjective, it has to be a covering map. If $G$ is abelian, $yz:=\alpha(\alpha^{-1}(y)\alpha^{-1}(z))$
is well-defined and turns $X\setminus D_{\textswab g}$ itself into an abelian Lie group of dimension $n$ and $\alpha$
into a group homomorphism. Since elements of $\ker\alpha$ induce the same action on $X\setminus D_{\textswab g}$,
the property $G\subset \Aut^0(X)$ implies, that $\alpha$ is an isomorphism. Hence we will identify
$G$ and $X\setminus D_{\textswab g}$ from now on.

Note that for a Lie group $G$ being abelian implies that $G={\bb C}^n/\Lambda$, where $\Lambda$ is a
discrete subgroup. This is proved by looking at the exponential map
$$\exp:\textswab g\seq G,$$
which is easily seen to be a group homomorphism of $(\textswab g,+)$ into $G$. Since $\exp$ maps some
neighborhood of $0$ diffeomorphically to a neighborhood of $1$, say $U$, and $\bigcup_{k=1}^\infty U^k=G$,
the map $\exp$ is surjective. Hence $G=(\textswab g,+)/\ker(\exp)={\bb C}^n/\Lambda,$ where $\Lambda:=\ker(\exp)$
must be discrete, since $n=\dim G=\dim\textswab g$.

In particular we see that $X\setminus D_{\textswab g}={\bb C}^n/\Lambda$, if $G$ is abelian. We will use this in the
next proof.

\begin{lemma}\label{bconstr}
Let $X$ be a compact complex manifold and $G\subset \Aut^0(X)$ a divisorial Lie group. 
If $g$ is a K\"ahler metric on $X\setminus D_{\textswab g}$ such that $\Aut^0(X,g)=G$
then $g$ is complete and there is a basis ${\mathcal B}\subset\textswab g$ such that 
$g=g_{\mathcal B}$.
\end{lemma}

\begin{proof}We already know by Lemma \ref{abel} that $G$ is 
abelian and hence $g_{\mathcal B}$ is K\"ahler and $G$-invariant,
if ${\mathcal B}\subset\textswab g$ is a basis. Since $G=X\setminus D_{\textswab g}={\bb C}^n/\Lambda$ we choose the images
of the canonical coordinates $z_1,...,z_n$ of ${\bb C}^n$ as local coordinates of $X\setminus D_{\textswab g}$. For
the sake of simplicity we call them also $z_1,...,z_n$. 
Of course, $g=g_{\alpha\bar\beta}dz^\alpha\otimes d\bar{z}^\beta$ is $G$-invariant,
if and only if $g_{\alpha\bar\beta}$ is constant for all $\alpha,\beta$. Hence $g$ is complete and 
corresponds
one to one to $g(0)$ what we identify with the matrix $\underline g=(g_{\alpha\bar\beta}(0))$.
The corresponding matrix $\underline g_{\mathcal B}$ is $\underline g_{\mathcal B}=\sigma\sigma^*$. Note that
$\sigma$ is constant since $G$ is abelian (cf. proof of Theorem \ref{binding}). Recall $S=\sigma^{-1}$
and define $H:=S\underline gS^*.$
Since $H$ is hermitian, we can find $A\in Gl(n)$ such that $H=A^*A$. Now 
$\underline g=\sigma H\sigma^*=\sigma AA^*\sigma^*.$
Set $B:=A^{-1}$. Then $g$ is given by the vector fields $t^{(i)}=\sum_jb_{ij}s^{(j)}$, 
which form another basis
of $\textswab g$. (Indeed, this shows by Theorem \ref{complete} once more that $g$ is complete.)
\end{proof}

\begin{cor}Let $X$ be a compact complex manifold and $S\subset X$ analytic with $\codim S>1$. If $X$ allows
for a divisorially invariant K\"ahler metric on $X\setminus S$, then $X$
is a torus.
\end{cor}

\begin{proof}Assume $g$ is such a metric and $G$ the divisorial abelian Lie group. By Lemma \ref{bconstr} $g|X\setminus D_{\textswab g}$ is constructed by a basis of $\textswab g$.
If $D_{\textswab g}$ is given by $\sigma\in H^0(-K_X)$, then $\det g=|\sigma|^{-2}$, hence is singular on $D_{\textswab g}$. 
This implies $D_{\textswab g}=0$. In particular, $X=G={\bb C}^n/\Lambda$. Since $X$ is compact,
$\Lambda$ is a complete lattice and $X$ is a torus.
\end{proof}
 
Note that this proof works also, if $\codim S=1$, but $S\not=D_{\textswab g,red}$. 

It is generally known (earliest references refer to a lecture of Remmert in 1958/1959) 
that any abelian complex Lie group is a direct product of copies of 
$T, {\bb C}$ and ${\bb C}^*$, where $T$ is a group without non-constant holomorphic functions.  
For a more detailed analysis we refer to \cite{mo}. Since Winkelmann related the triviality of a logarithmic tangent bundle to $G$
being a semi-torus in \cite{wi}, we would like to express this property in terms of the mentioned decomposition. 

\begin{lemma}\label{st}$G=T\times({\bb C}^*)^k\times{\bb C}^l$ is a semi-torus $\iff l=0$.\end{lemma}

\begin{proof}If $G$ is not a semi-torus, i.e. $\Lambda\otimes{\bb C}\not={\bb C}^n$, then obviously $l\not=0$. So let $G$ be a semi-torus. 
There is a lattice $\Lambda'\subset{\bb C}^{n-l}$ coming from the decomposition such that $G={\bb C}^n/\Lambda'$. The isomorphism 
$\tilde\phi:{\bb C}^n/\Lambda\seq{\bb C}^n/\Lambda'$ is induced by a vector space isomorphism $\phi:{\bb C}^n\seq{\bb C}^n$ obeying $\phi(\Lambda)=\Lambda'$.
Hence ${\bb C}^{n-l}=\Lambda'\otimes{\bb C}=\phi(\Lambda\otimes{\bb C})={\bb C}^n$, so $l=0$.
\end{proof}

The toric varieties fit in the system
as special cases where $G=({\bb C}^*)^n$. 
If $X$ is Fano, then $X\setminus D_{\textswab g}$ is Stein and hence
the factor $T$ does not occur (cf. also
\cite[Prop. 4]{mm} for this particular claim). 
If we relax this condition a little we can show that $T$ is a torus. The following lemma in particular covers the case of homogeneous manifolds for which the result
is well-known.

\begin{lemma}Let $X$ be a projective almost homogeneous manifold such that $|-mK_X|$ is base point free for a certain $m>0$. Let further
$G\subset \Aut^0(X)$ be an abelian divisorial Lie group.
$X\setminus D_{\textswab g}$ is of the form $T\times{\bb C}^k\times({\bb C}^*)^l$, where $T$ is a 
torus. Furthermore, $X=T\times Y$ for a rational manifold $Y$.
\end{lemma}

\begin{proof}We first have to prove that $T$ is a torus. The assumption that $|-mK_X|$ is base point free enables us to choose for every $x\in D_{\textswab g}$
a meromorphic function $\tilde f$ with poles exactly along $D_{\textswab g}$ and $x$ is not in the locus of indeterminacy of $f$. If $T$ is not compact, we
fix $z\in{\bb C}^k\times({\bb C}^*)^l$ and $x\in\overline{T\times\{z\}}$.
Since $\tilde f|_{X\setminus D_{\textswab g}}$ is holomorphic, 
in particular $f:=\tilde f|_{T\times\{z\}}$ is holomorphic. 
Hence 
$f$ is constant and we obtain $\overline{T\times\{z\}}\subset X\setminus D_{\textswab g}$. This implies that $T$ is compact. 

Hence $T$ is a projective manifold.
Since $T={\bb C}^k/\Lambda$, the lattice $\Lambda$ is complete and hence $T$ is a torus.

Now we have to prove that the projection onto $T$ is extendable. Since $X$ is bimeromorphic to
${\bb P}^n\times T$, and the Albanese torus $A(X)$ as well as the Albanese map $\alpha:X\seq A(X)$ are
bimeromorphic invariants of projective manifolds, $T=A(X)$ and $pr=\alpha$ by the universal property of $(A(X),\alpha)$.

Now choose $h:=(1,t')\in G=H\times T$ and denote $F_t:=\alpha^{-1}(t)$. Of course, $h:X\seq X$ satisfies
$h(F_t)=F_{t+t'}$ and the map $\psi:F_0\times T\seq X, (y,t)\mapsto (1,t)y$ is an isomorphism. Since $Y:=F_0$ is a
compactification of ${\bb C}^k\times({\bb C}^*)^l$, it is rational.
\end{proof}

Indeed, all factors of $G$ can occur. If $X$ factors $X=T\times\tilde X$, then 
$X\setminus D_{\textswab g}=T\times\tilde X\setminus\tilde D$, since $K_T=0$. The factor ${\bb C}^k\times({\bb C}^*)^l$
occurs even for $X={\bb P}^n$, if $G$ is carefully chosen. For example, the group
$$G=\left\{\left.\left(\begin{array}{ccc}1 & t & s+\frac{t^2}{2}\\
0& 1 & t \\0 & 0 & 1\end{array}\right)\right| s,t\in{\bb C}\right\}\cong({\bb C}^2,+)$$
acts on ${\bb P}^2$, leaving only $D_{\textswab g,red}=\{z^2=0\}$ invariant. Of course, the corresponding vector fields
$s^{(1)}=z^2\frac{\partial}{\partial z^0}, s^{(2)}=z^2\frac{\partial}{\partial z^1}+
z^1\frac{\partial}{\partial z^0}$ yield $D_{\textswab g}=\{(z^2)^3=0\}$. 
We leave it as an exercise to the reader to construct
the other cases. The philosophy is: ${\bb C}^*$-actions degenerate to ${\bb C}$-actions whenever two 
hypersurfaces coincide. Indeed, we are now going to show that this is always the situation.

\begin{lemma}If $D_{\textswab g}$ is reduced, then $G$ is a semi-torus.
\end{lemma}

\begin{proof}We assume that $G$ is not a semi-torus, hence by Lemma \ref{st} we obtain a decomposition $G={\bb C}\times G'$. 
We choose $s\in H^0(-K_X)$ such that $D_{\textswab g}=\{s=0\}$. Since $G'$ is abelian,
$G'={\bb C}^{n-1}/\Lambda$ and we may choose local coordinates $\tilde z^1,...,\tilde z^{n-1}$ induced
by canonical coordinates of ${\bb C}^{n-1}$. In those coordinates of $X\setminus D_{\textswab g}$
we write
$$s=f(x,\tilde z^1,...,\tilde z^{n-1})\frac{\partial}{\partial x}\wedge\frac{\partial}{\partial \tilde z^1}
\wedge ...\wedge\frac{\partial}{\partial\tilde z^{n-1}},$$
where $x$ denotes the coordinate of the factor ${\bb C}$. If $x\seq\infty$, we will approximate
a point in $D_{\textswab g}$. In order to approximate other points (and indeed by this procedure
all other points of a certain
component of $D_{\textswab g}$), we choose an arbitrary holomorphic $\lambda:{\bb C}\seq G'$ and look at the curve
$(x,\tilde z'+\lambda(x))$ for a fixed point 
$\tilde z':=(\tilde z^1,...,\tilde z^{n-1})\in G'$. Let $p=\lim_{x\seq\infty}(x,\tilde z'+\lambda(x))\in D_{\textswab g}$.
By $x\mapsto\frac 1x=:y, \tilde z^i\mapsto\tilde z^i-\lambda^i(x)=:z^i$
we get local coordinates
in a punctured neighborhood $U(p)\setminus\{p\}$. In these coordinates,
$$s=-f\left(\frac 1y,z'+\lambda\left(\frac 1y\right)\right)y^2\frac{\partial}{\partial y}\wedge\frac{\partial}{\partial z^1}\wedge ...
\wedge\frac{\partial}{\partial z^{n-1}}.$$ Let us denote $h(y):=-f(\frac 1y,z'+\lambda(\frac 1y))y^2.$ 
The group action
of ${\bb C}$ now is 
$$\lambda\cdot y=\frac{y}{1+\lambda y}.$$
The invariance of $D_{\textswab g}$ under $G$ implies for $\lambda\in{\bb C}$ that
$\lambda^*s=c(\lambda)s$, hence
$$h(\lambda\cdot y)(1+\lambda y)^2=c(\lambda)h(y).$$
Since $c(\lambda+\kappa)s=(\lambda+\kappa)^*s=\lambda^*\kappa^*s=c(\lambda)c(\kappa)s,$
the function $c(\lambda)=\exp(\rho\lambda).$ This implies
$$h(\lambda\cdot y)=\frac{\exp(\rho\lambda)}{(1+\lambda y)^2}h(y).$$
Now fixing $y=1$ yields
$$h(\frac 1{1+\lambda})=c\frac{\exp(\rho(1+\lambda))}{(1+\lambda)^2},$$
hence
$$h(y)=cy^2\exp(\frac{\rho}{y}).$$
Since we have the additional requirement that $s|D_{\textswab g}=0$ and $\frac{\partial}{\partial y}\wedge\frac{\partial}{\partial z^1}\wedge ...
\wedge\frac{\partial}{\partial z^{n-1}}(p)$ has a finite vanishing order, we conclude
$\rho=0$. Now we see that $h$ vanishes of order $2$ in $0$. We cannot guarantee that different choices
of $z'$ lead to different limit points on $D_{\textswab g}$, i.e. maybe 
$\frac{\partial}{\partial y}\wedge\frac{\partial}{\partial z^1}\wedge ...
\wedge\frac{\partial}{\partial z^{n-1}}(p)=0$.
Hence we only conclude that the vanishing order of $s$ on
the limit point $p$ is at least $2$. Since we could do this construction for every point of a
component containing $p$, we conclude that this component is multiple. 
\end{proof}

 
\begin{cor}\label{toric}
If $X$ is an $n$-dimensional projective compact complex manifold such that $|-mK_X|$ is base point free for some $m>0$, further $D\in|-K_X|$ a reduced divisor and
$g$ a K\"ahler metric on $X\setminus D$ with divisorial invariance group $\Aut(X,g)$,
then $X=T\times P$, where $P$ denotes a projective toric variety; in this case $D=\sum T\times D_i$, where
$D_i$ are the distinct toric divisors of $P$.
\end{cor}

\begin{proof}We know $X=P\times T$. Since $D$ is reduced, $X\setminus D=({\bb C}^*)^l\times T$. 
Since $P$ is algebraic and has an algebraic $({\bb C}^*)^l$-action, it is a toric variety and $D$
is like described (cf. \cite{fu}).
\end{proof} 

This splitting behaviour cannot be expected in general. However, $X$ is always a fibre bundle over $\Alb(X)$. 
This is generally known and easy to be seen by the universal property of $\alpha$.

\subsection{Non-Triviality of the $G$-K\"ahler cone}

Recall the definition of the $G$-K\"ahler cone: We call two $G$-invariant
K\"ahler metrics $g$ and $g'$ 
equivalent, if there
is a function $\phi\in C^\infty(X\setminus D_{\textswab g})$ such that $\omega_g-\omega_{g'}=i\dd\phi$. Since there is 
no $\dd$-lemma
in the non-compact case, this cannot be viewed as the K\"ahler class of $\omega_g$, but the philosophy 
is very similar. Hence we denote $K_G(X,D_{\textswab g}):=M_G(X,D_{\textswab g})/\sim$.

We identify again $X\setminus D_{\textswab g}={\bb C}^n/\Lambda$. Every $\lambda_{(1)},\lambda_{(2)}\in\Lambda$ generate a
parallelogram, whose image in $X\setminus D_{\textswab g}$ is a compact real surface $T_{\lambda_{(1)},\lambda_{(2)}}$.
If $\omega=i\dd\phi$, then Stokes' Theorem implies
$$0=\int_{T_{\lambda_{(1)},\lambda_{(2)}}}\omega=\omega_{ij}(\lambda_{(1)}^i\overline{\lambda_{(2)}^j}-
\lambda_{(2)}^i\overline{\lambda_{(1)}^j})=2\im(\lambda_{(1)}^t\underline{\omega}\overline{\lambda_{(2)}}),$$
if $\underline\omega=(\omega_{ij})_{i,j}$. Note that the $\omega_{ij}$ are constant.
This property only depends on $\Lambda_{\bb R}:=\Lambda\otimes{\bb R}$. It is easy to see that in appropriate complex coordinates every
real subspace of ${\bb C}^n$ is of the form
$$\Lambda_{\bb R}=\{z^1=\ldots=z^{l'}=\im z^{l'+1}=\ldots=\im z^{k'}=0\}.$$
In other words, $\Lambda_{\bb R}$ is generated
by the real standard basis of ${\bb C}^k={\bb R}^{2k}$ and the standard basis of ${\bb R}^l=\re({\bb C}^l)$ (for $k=n-k', l=k'-l'$).
By this choice the above equations mean that in the standard basis of ${\bb C}^n={\bb C}^k\oplus{\bb C}^l\oplus{\bb C}^{n-k-l}$
$$\underline\omega=\left(\begin{array}{ccc}0 & 0 & *\\0 & \mbox{real} & *\\ * & * & *\end{array}\right),$$ 
where every entry stands for the block corresponding to the factors of ${\bb C}^n={\bb C}^k\oplus{\bb C}^l\oplus{\bb C}^{n-k-l}$ and $*$ means,
that there is no claim about this entry. 

Now let us reverse the direction.
For the sake of simplicity, let us denote
${\bf k}:=\{1,...,k\}$, ${\bf l}:=\{k+1,...,k+l\}$, ${\bf m}:=\{k+l+1,...,n\}$.
If $\omega$ is of the above form, we define $\phi$ by
\begin{eqnarray*}\phi(z)&:=&2\sum_{i\in\bf l}\omega_{ii}\im(z^i)^2+4\sum_{i<j\in\bf l}\omega_{ij}\im(z^i)\im(z^j)+\\
& & +4\sum_{i\in {\bf l},j\in {\bf m}}(\omega_{ij}\im(z^i)\overline{z^j}-\omega_{ji}\im(z^i)z^j)+\\
& & +4\sum_{i\in {\bf m}}\omega_{ii}|z^i|^2
+2\sum_{i<j\in {\bf m}}\omega_{ij}z^i\overline{z^j},\end{eqnarray*}
then
$$\omega-i\dd\phi=\left(\begin{array}{ccc}0 & 0 & *\\0 & 0 & 0\\ * & 0 & 0\end{array}\right).$$

If we now assume further, that $G$ is semi-torus (e.g. $D_{\textswab g}$ is reduced), then the factor ${\bb C}^{n-k-l}$ does not occur and hence
$$\omega=i\dd\phi.$$

This leads to the result

\begin{thm}\label{cone}Let $X$ be an $n$-dimensional compact complex manifold and $G={\bb C}^n/\Lambda$ a divisorial semi-torus. 
We write $\Lambda_{\bb R}={\bb C}^k\oplus\re({\bb C}^l)$ (with $k+l=n$).
Denote $i:M(l,{\bb C})\seq M(n,{\bb C})$ the embedding which fills up
an $l\times l$-matrix with zeroes. Then   
$$K_G(X,D_{\textswab g})\subset M(n,{\bb C})/i(M(l,{\bb R}))$$
is the cone generated by positive definite hermitian matrices. In particular, 
$\dim K_G(X,D_{\textswab g})=n^2-\frac 12l(l+1).$
\end{thm}

This contrasts to the case of a smooth divisor $D$. In the appendix we will show that the K\"ahler cone is
trivial, if $X$ has simple topology and $D$ is smooth.

\section{Example: $X={\bb P}^2$}\label{example}

If $X={\bb P}^2$, then the tangent bundle may be described by the vector fields homogeneous of degree $1$ divided
by the vector fields parallel to the orbits of the group action $z\mapsto cz$, i.e. 
${\mathcal O}_X\cdot(z^0\frac{\partial}{\partial z^0}+z^1\frac{\partial}{\partial z^1}+z^2\frac{\partial}{\partial z^2})$.
Hence the global vector fields are
$$H^0(T_X)\cong (l^0\frac{\partial}{\partial z^0}+l^1\frac{\partial}{\partial z^1}+
l^2\frac{\partial}{\partial z^2})/{\bb C}\cdot(z^0\frac{\partial}{\partial z^0}+z^1\frac{\partial}{\partial z^1}+
z^2\frac{\partial}{\partial z^2}),$$
where $l^i$ are homogeneous linear forms.
Now let $V:={\bb C}v^{(1)}\oplus{\bb C}v^{(2)}\subset H^0(T_X)$ with 
$v_j=[\sum_il^{ji}\frac{\partial}{\partial z^i}]$. 
In order to compute
$$D:=\{z|v^{(1)}\wedge v^{(2)}=0\},$$
we first localize to $U_0$ and then homogenize the result again. This procedure yields
$$D=\left\{\det\left(\begin{array}{ccc}z^0&z^1&z^2\\l^{10}&l^{11}&l^{12}\\l^{20}&l^{21}&l^{22}\end{array}\right)=0\right\}.$$

Now let us assume that $[v^{(1)},v^{(2)}]=0$ and $v^{(1)}\wedge v^{(2)}\not\equiv  0$. Denote $G:=\exp(V)=
\Aut^0(X,D)$. By assumption $G$ is divisorial and abelian, 
hence the metric $g_{\mathcal B}$ is K\"ahler and $G=\Aut(X,g_{\mathcal B})$. 
If $D$ is reduced, Theorem \ref{toric} tells us that $D$ is the 
union of three lines in general position. If $D$ is not reduced we obtain two lines one of which is double or
a triple line. So the only position of three lines not occurring in this list is that they are intersecting
in one common point. We will now see how this corresponds to a $G$ which acts not almost transitively.
After a change of coordinates we may assume that the three lines intersect in $[1:0:0]$.
Let $v^{(1)}:=z^1\frac{\partial}{\partial z^0}, v^{(2)}:=z^2\frac{\partial}{\partial z^0}$. Of course,
$[v^{(1)},v^{(2)}]=0$ and
$v^{(1)}\wedge v^{(2)}\equiv 0$, hence $G$ is abelian (indeed, $G\cong{\bb C}^2$) and
acts not almost transitively. $G$ is given by the matrices
$$\left(\begin{array}{ccc}1 & a & b\\0 & 1 & 0\\0 & 0 & 1\end{array}\right), a,b\in{\bb C}.$$
It is not hard to see that $G$ leaves $\{f=0\}$ invariant for a homogeneous
$f\in{\bb C}[z^0,z^1,z^2]$ if and only if $f=f(z^1,z^2)$. This factors into linear terms. Hence
$G$ leaves all lines through $[1:0:0]$ invariant. So the not almost transitively case corresponds
to the existence of a family of invariant divisors, which are not necessarily anticanonical. 

Back to the almost transitive $G$ and reduced $D$. Let us choose coordinates such that
$D=\{z^0z^1z^2=0\}$. Of course, $D$ is invariant under the group $G$ given by $[z^0:z^1:z^2]\mapsto
[a_0z^0:a_1z^1:a_2z^2]$, with $a=[a_0:a_1:a_2]\in{\bb P}^2\setminus\{a_0a_1a_2=0\}
\cong{\bb C}^*\times{\bb C}^*$. The group $G$ is abelian and divisorial.
Theorem \ref{binding} and Lemma \ref{prebinding} tell us that $G=\Aut^0(X,D)$. Lemma \ref{bconstr}
states that every $G$-invariant K\"ahler metric on $X\setminus D$ is given by a basis of
$\textswab g$. Carrying out the calculations in the chart $U_0=\{z_0\not=0\}$ yields that all $G$-invariant
K\"ahler metrics on $X\setminus D\cong{\bb C}^*\times{\bb C}^*$ are of the form
$$g=g_C=\sum_{i,j=1,2} c_{ij}\frac{dz^i}{z^i}\otimes\frac{d\overline{z^j}}{\overline{z^j}},$$
with $c_{ij}=\overline{c_{ji}}$ and $C=(c_{ij})>0$.

According to Theorem \ref{cone}
$$K_G(X,D)\subset M(2,{\bb C})/M(2,{\bb R})$$
is given by the classes  of positive, hermitian matrices. Hence $K_G(X,D)$ is one-dimensional. Moreover
$$C(r):=\left(\begin{array}{cc}\cosh(r)&i\sinh(r)\\-i\sinh(r)&\cosh(r)\end{array}\right)$$
for $r\in{\bb R}$ represent every class in $K_G(X,D)$ uniquely.

If $X={\bb P}^3$, then for the corresponding construction $\dim K_G({\bb P}^3,D)=3$. 
In the appendix it is proved that $K({\bb P}^3,D)=0$, if $D$ is chosen to be smooth.

\appendix\section{Triviality of the K\"ahler cone when $D$ is smooth}

\begin{thm}\label{trivial}
Let $X$ be a projective complex manifold with $\dim X\ge 3, b_1(X)=b_3(X)=0, b_2(X)=1$ and
$D\subset X$ a smooth ample divisor. Then $K(X\setminus D)=0$.
\end{thm}

The proof will be divided into three steps. We abbreviate $\tilde X:=X\setminus D$. Of course,
since $\tilde X$ is quasi-projective, $\tilde X$ is K\"ahler. So $K(X,D)\not=\emptyset$.

\begin{lemma}\label{1}Under the assumptions of Theorem \ref{trivial} holds $b_1(\tilde X)=0$.\end{lemma}

\begin{proof}
Since $D$ is ample, $\tilde X$ is Stein. Since ${\mathcal O}_{\tilde X}$ is coherent, Theorem B
implies $H^1(\tilde X,{\mathcal O})=0$. Now let $\varphi\in{\mathcal E}^1(\tilde X)$ be a closed
real one-form. If we decompose into types $\varphi=\varphi^{0,1}+\varphi^{1,0}$, then we obtain
$$\partial\varphi^{1,0}=\overline\partial\varphi^{0,1}=
\overline\partial\varphi^{1,0}+\partial\varphi^{0,1}=0.$$
Dolbeault cohomology yields functions $g,h:\tilde X\seq{\bb C}$ such that $\partial g=\varphi^{1,0},
\overline\partial h=\varphi^{0,1}$. Now we compute
$$d(g+h)=\varphi+\overline\partial g+\partial h.$$
Since
$$\partial\overline\partial g=-\overline\partial\varphi^{1,0}=\partial\varphi^{0,1}=\partial\overline\partial h,$$
we conclude that there exists a holomorphic function $F\in H^0({\mathcal O}_{\tilde X})$ such that
$$g-h=F+\overline F.$$
This allows the computations
$$\overline\partial g+\partial h=\partial h+\overline\partial(h+F+\overline F)=\partial h+\overline\partial
\,\overline F+\varphi^{0,1}$$
and
$$\overline\partial g+\partial h=\overline\partial g+\partial(g-F-\overline F)=\overline\partial g+
\varphi^{1,0} -\partial F.$$
Adding both equations yields
$$\overline\partial g+\partial h=\varphi+d(\overline F-F),$$
hence
$$d(g+h)=2\varphi+d(\overline F-F),$$
hence $\varphi$ is $d$-exact.
\end{proof}

\begin{lemma}\label{2}Under the assumptions of Theorem \ref{trivial} holds $b_2(\tilde X)=0$.\end{lemma}

\begin{proof}
We choose a tubular neighborhood of $D$ and use Mayer-Vietoris for $X=\tilde X\cup U(D)$.
Since $\tilde X\cap U(D)=U(D)\setminus D$ has a $S^1$-bundle $E\seq D$ as a deformation retract
and $U(D)$ is contractible to $D$, Mayer-Vietoris yields
\begin{eqnarray*}H^1(X,{\bb C})\seq H^1(\tilde X,{\bb C})\oplus H^1(D,{\bb C})\seq H^1(E,{\bb C})
\seq H^2(X,{\bb C})& &\\ 
\seq H^2(\tilde X,{\bb C})\oplus H^2(D,{\bb C})
\seq H^2(E,{\bb C})\seq H^3(X,{\bb C}).& &\end{eqnarray*}
In order to compute the cohomology of $E$ we use the K\"unneth formula for locally trivial
fibrations (cf. \cite[p. 258]{spa}). This yields
$$b_1(E)=b_1(D)+b_0(D).$$
Our assumptions that $D$ is smooth and ample and $\dim X\ge 3$ allow us to use the Lefschetz theorem
to conclude $h^1(D,{\bb C})=0$. Hence $b_1(E)=1$.
In the same way we compute $b_2(E)=b_2(D)$.
If we now use the assumptions $b_1(X)=b_3(X)=0, b_2(X)=1$ as well as $b_1(\tilde X)=0$ by
Lemma \ref{1}, then the sequence implies $b_2(\tilde X)=0$.
\end{proof}

Now we proceed to the proof of the theorem.

\begin{proof}[of Theorem \ref{trivial}]
The injective resolution of ${\bb C}$
$$0\seq{\bb C}\seq{\mathcal O}_{\tilde X}\stackrel{\partial}{\seq}\Omega^1_{\tilde X}\stackrel{\partial}{\seq}
\Omega^2_{\tilde X}\stackrel{\partial}{\seq}...$$
yields short exact sequences
$$0\seq{\bb C}\seq{\mathcal O}_{\tilde X}\stackrel{\partial}{\seq}{\mathcal H}^1_{\tilde X}\seq 0$$
and
$$0\seq{\mathcal H}^1_{\tilde X}\seq\Omega^1_{\tilde X}\seq{\mathcal H}^2_{\tilde X}\seq 0.$$
In cohomology we obtain
$$H^1({\mathcal O}_{\tilde X})\seq H^1({\mathcal H}^1_{\tilde X})\seq H^2({\tilde X},{\bb C})\seq 
H^2({\mathcal O}_{\tilde X}).$$
Since $\tilde X$ is Stein we obtain $H^1({\mathcal O}_{\tilde X})=H^2({\mathcal O}_{\tilde X})=0,$
hence
$$H^1({\mathcal H}^1_{\tilde X})=H^2(\tilde X,{\bb C})=0$$
by Lemma \ref{2}.
The second short exact sequence yields
$$H^0(\Omega^1_{\tilde X})\seq H^0({\mathcal H}^2_{\tilde X})\seq H^1({\mathcal H}^1_{\tilde X})=0,$$
hence for every holomorphic 
$2$-form $\eta$ on $\tilde X$ with $\partial\eta=0$ there is a holomorphic $1$-form
$\varphi$ such that $\eta=\partial\varphi$.

Now let $\omega$ be a K\"ahler form. Since $\Omega^1_{\tilde X}$ is coherent, again Theorem B implies
$H^1(\Omega^1_{\tilde X})=0$. Using the Dolbeault interpretation we obtain
$\eta\in{\mathcal E}^{1,0}(\tilde X)$ such that $\omega=\overline\partial\eta$.
Now look at $\psi:=\partial\eta$. Since $0=\partial\omega=-\overline\partial\psi$, we conclude that
$\psi\in H^0(\Omega^2_{\tilde X}).$  Of course, it satisfies $\partial\psi=0$.
Hence there is a $\varphi\in H^0(\Omega^1_{\tilde X})$ such that
$\psi=\partial\varphi$. This implies $\partial(\eta-\varphi)=0$, hence $\overline\eta-\overline\varphi$
induces a class in $H^{0,1}(\tilde X)=H^1({\mathcal O}_{\tilde X})=0$. Hence we obtain a function
$G:\tilde X\seq{\bb C}$ such that $\overline\partial G=\overline\eta-\overline\varphi$, hence
$$\partial\overline G=\eta-\varphi.$$
For the K\"ahler form this means
$$\omega=\overline\partial\eta=\overline\partial(\partial\overline G+\varphi)=\partial\overline\partial
(-\overline G).$$
Since $\omega=\overline{\omega}$ we find
$$\omega=i\partial\overline\partial\im(G).$$
\end{proof}

\

\noindent
{\bf Acknowledgments.} The authors are grateful to Hans-Christoph Grunau and Georg Schumacher 
for proposing and supporting the project.

\end{document}